\documentclass[10pt,reqno]{amsart}

\usepackage[T1]{fontenc}
\usepackage{lmodern}
\usepackage{color,xcolor}
\usepackage{amsmath,amssymb,amsthm}
\usepackage[letterpaper,margin=1in]{geometry}
\usepackage{microtype}
\usepackage[backend=biber,style=numeric]{biblatex}
\usepackage[colorlinks=true,linkcolor=blue,citecolor=blue,urlcolor=blue]{hyperref}

\newtheorem{theorem}{Theorem}
\newtheorem{lemma}[theorem]{Lemma}

\theoremstyle{definition}

\theoremstyle{remark}

\newcommand{\R}{\mathbb R}
\newcommand{\N}{\mathbb N}
\newcommand{\Nzero}{\mathbb N_0}
\newcommand{\norm}[1]{\lVert #1\rVert}
\newcommand{\dual}[2]{\langle #1,#2\rangle}

\title[Quasi-orthogonality
is a geometric property]{Quasi-orthogonality is a geometric property}
\date{}
\hypersetup{
  pdftitle={Quasi-orthogonality for strongly monotone Galerkin methods},
  pdfkeywords={strong monotonicity, nonlinear Galerkin methods, quasi-orthogonality}
}
\author{Michael Feischl and David Niederkofler}
\thanks{The research for this work was heavily AI supported and we provide a LEAN certificate of Theorem~\ref{thm:qo} under \url{https://github.com/michaelFeischl/quasi-orthogonality-lean}.
Funded by the Deutsche Forschungsgemeinschaft (DFG, German Research Foundation) -- Project-ID 258734477 -- SFB 1173, the Austrian Science Fund (FWF)
under the special research program Taming complexity in PDE systems (grant SFB F65) as well as projects I6667-N and P36150. This work was also supported by the ERC grant OPTIMAL (grant agreement No. 101125225). }
\begin{document}
\begin{abstract}
    Quasi-orthogonality is one of the key ingredients for proofs of optimal convergence rates for adaptive mesh refinement methods~\cite{axioms}. The work~\cite{infsup} proposed a relaxed quasi-orthogonality that follows directly from the inf-sup stability of the underlying problem discretization. The proof uses an intricate connection to the stability of $LU$-factorization of matrices. Recently,~\cite{gantumur} presented a vastly simplified argument based solely on Hilbert space arguments. Building on this idea, we show that quasi-orthogonality is a purely geometric property of a set of points. This property is checked for a large class of Petrov-Galerkin methods for monotone, non-linear equations.
\end{abstract}
\maketitle

\section{Main result}
We consider a real Hilbert space $U$ and a finite sequence of elements $u_1,\ldots,u_n\in U$. For each subset $J\subseteq \{1,\ldots,n\}$ and $0\leq \lambda_j\leq 1$ with $\sum_{j\in J}\lambda_j=1$ we define the weighted mean
\begin{align*}
    \mu_{J,\lambda}:= \sum_{j\in J}\lambda_j u_j
\end{align*}
and the maximal weighted squared distances (variance)
\begin{align*}
     E_J : = \max_{0\leq \lambda_i\leq 1\atop \sum_{j\in J} \lambda_j = 1}\sum_{j\in J}\lambda_j \norm{u_j- \mu_{J,\lambda}}_U^2.
\end{align*}
We make the following two assumptions:
\begin{itemize}
    \item C\'ea type estimate: There exists $C_{\text{C\'ea}}>0$ such that for all $1\leq i\leq j\leq k\leq \ell\leq n$, there holds
    \begin{align}\label{eq:cea}
    \norm{u_k-u_j}_U\leq C_{\text{C\'ea}}\norm{u_\ell-u_i}_U.
    \end{align}
\item Separation: There exists $C_{\rm sep}>0$ such that for all $J\subseteq \{1,\ldots,n\}$, all $v\in {\rm conv}\{u_i\,|\,1\leq i\leq n,\,i\leq J\}$, there holds
\begin{align}\label{eq:sep}
        C_{\rm sep}\norm{\mu_{J,\lambda}-v}_U^2 \geq   E_J,
    \end{align}
    where $(\lambda_j)_{j\in J}$ is the maximizer of $E_J$.
\end{itemize}

For such sequences, we can prove the following quasi-orthogonality result.
\begin{theorem}\label{thm:qo}
There exists $C>0$ that depends only on the constants $C_{\text{C\'ea}}$, and $C_{\rm sep}$ such that
\begin{align*}
    \sum_{i=1}^{n-1} \norm{u_{i+1}-u_i}^2_U\leq C n^{1-\delta} \norm{u_n-u_1}_U^2,
\end{align*}
where $\delta := 1-\log_2\max\{\frac32,(\frac12 + \frac{1}{16C_{\rm sep}})^{-1}\}>0$.
\end{theorem}
The main idea of the following proof is to look at the variance (mean squared distance from the center) of different subsets of $u_1,\ldots, u_n$.
By using the fact that $u_{i+1},\ldots,u_n$ is geometrically separated from $u_1,\ldots, u_i$ for all $1\leq i\leq n$, we will show that combining two disjoint subsets increases the variance by a factor larger than $1/2$ over the sum of the individual variances. In Galerkin methods, this separation can often be established by the fact that the residual of the equation w.r.t. $u_j$, $j\geq i$ is orthogonal to $u_1,\ldots, u_i$ (see Section~\ref{sec:pg} for details).
Moreover, the variance of subsets of the form $\{i,i+1\}$ is an upper bound for $\norm{u_{i+1}-u_i}_U^2$. Combining both observations allows us to iteratively combine two disjoint sets into larger ones until we reach the full set. This requires roughly $\log_2(n)$ iterations. Since the right-hand side bound is the inverse product of the growth factors in each iteration, we obtain the result. 

\begin{proof}[Proof of Theorem~\ref{thm:qo}]
    \textbf{Step 1:} 
    There holds
    \begin{align}
        E_I = \max_{0\leq \lambda_i\leq 1\atop \sum_{i\in I} \lambda_i = 1}\Big(\sum_{i\in I}\lambda_i \norm{u_i}_U^2  - 2\sum_{i\in I}\lambda_i \dual{u_i}{\mu_{I,\lambda}}_U+\norm{\mu_{I,\lambda}}_U^2\Big) = \max_{0\leq \lambda_i\leq 1\atop \sum_{i\in I} \lambda_i = 1} \Big(\sum_{i\in I}\lambda_i \norm{u_i}_U^2  - \norm{\mu_{I,\lambda}}_U^2\Big). 
    \end{align}
    Given $I,J\subseteq \{1,\ldots,n\}$ with $I\cap J = \emptyset$ and $i < j$ for all $i \in I$ and $j \in J$, with corresponding maximizers $\lambda_I$ and $\lambda_{J}$, there holds
    \begin{align}\label{eq:Eest}
    \begin{split}
        E_{I\cup J} &\geq  \frac12\Big(\sum_{i\in I}\lambda_{I,i} \norm{u_i}_U^2+\sum_{i\in J}\lambda_{J,i} \norm{u_i}_U^2\Big)-\frac14\norm{\sum_{i\in I}\lambda_{I,i}u_i+\sum_{i\in J}\lambda_{J,i}u_i}_U^2\\
        &\geq  \frac12\Big(E_I + E_{J}\Big)+\frac14\norm{\sum_{i\in I}\lambda_{I,i}u_i}_U^2+\frac14\norm{\sum_{i\in J}\lambda_{J,i}u_i}_U^2
        -\frac12 \dual{\sum_{i\in I}\lambda_{I,i}u_i}{\sum_{i\in J}\lambda_{J,i}u_i}_U\\
        &=
        \frac12\Big(E_I + E_{J}\Big)+\frac14\norm{\mu_{I,\lambda_I}-\mu_{J,\lambda_{J}}}_U^2.
        \end{split}
    \end{align}

    \textbf{Step 2:}    
    Separation~\eqref{eq:sep} together with~\eqref{eq:Eest} shows for $c = 1/(4C_{\rm sep})$
    \begin{align}\label{eq:quadest}
        E_{I\cup J}\geq \frac12\Big(E_I + E_{J}\Big)+cE_{J}.
    \end{align}
    If $E_{J}\geq \frac12 E_I$, the above inequality~\eqref{eq:quadest} shows
    \begin{align*}
        E_{I\cup J}\geq (\frac{1}{2}+\frac{c}{4})(E_I+E_{J}).
    \end{align*}
    If $E_{J}<\frac12 E_I$, we use $E_I \leq E_{I \cup J}$ to also get $(E_I+E_{J})\leq \frac32 E_I\leq \frac32 E_{I\cup J}$.
    Thus, in both cases we have 
    \begin{align}\label{eq:reduction}
    E_I+E_{J}\leq qE_{I\cup J}
    \end{align}
    for $q:=\max\{3/2,(1/2+c/4)^{-1}\}<2$.

    \textbf{Step 3}: 
    There holds for $k=\lfloor\log_2(n-1)\rfloor$ that
    \begin{align*}
        \max\Big\{\sum_{i=1}^{2^k} \norm{u_{i+1}-u_i}^2_U, \sum_{i=n-2^k}^{n-1} \norm{u_{i+1}-u_i}^2_U\Big\}\geq \frac12 \sum_{i=1}^{n-1} \norm{u_{i+1}-u_i}^2_U.
    \end{align*}
    Thus, without loss of generality, we may assume $n-1=2^k$. Moreover, there holds
    \begin{align*}
        \max\Big\{\sum_{j=1}^{(n-1)/2} \norm{u_{2j}-u_{2j-1}}^2_U, \sum_{j=1}^{(n-1)/2} \norm{u_{2j+1}-u_{2j}}^2_U\Big\}\geq \frac12 \sum_{i=1}^{n-1} \norm{u_{i+1}-u_i}^2_U.
    \end{align*}
    Without loss of generality, we assume that the first sum in the maximum above is the larger one and start with $I_j:= \{2j-1,2j\}$ for all $j=1,\ldots,(n-1)/2$. Using~\eqref{eq:reduction}, we may iteratively combine consecutive $I_j$ and $I_{j+1}$ to obtain
    \begin{align*}
        \sum_{j=1}^{(n-1)/2}E_{I_j}\leq q^{\log_2(n)} E_{\{1,\ldots,n\}}.
    \end{align*}
    Clearly, we have $E_{\{1,\ldots,n\}}\lesssim  \norm{u_n-u_1}_U^2$ and also $E_{I_j}\geq \frac14\norm{u_{2j}-u_{2j-1}}_U^2$.
    This and $q^{\log_2(n)} = n^{\log_2(q)}$ conclude the proof.
\end{proof}

\section{Application to non-linear Petrov Galerkin Methods}\label{sec:pg}
\noindent

We show that the assumptions~\eqref{eq:cea} and~\eqref{eq:sep} are relatively straightforward for a large class of nonlinear Petrov-Galerkin methods.
Let $U,V$ be real Hilbert spaces. We denote the duality pairing between
$V^\star$ and $V$ by $\dual{\cdot}{\cdot}$, and equip each subspace occurring in the discussion
with the norm inherited from its ambient space. Suppose that
$A:U\to V^\star$ satisfies
\begin{equation}\label{eq:pg-lipschitz}
  \norm{Ax-Ay}_{V^\star}\le L\norm{x-y}_U,
  \qquad x,y\in U,
\end{equation}
for some $L>0$. Let $(U_j)_{j\ge0}$ and $(V_j)_{j\ge0}$ be nested
sequences of finite-dimensional linear subspaces of $U$ and $V$,
respectively. Nestedness means that for all $j\ge0$, we have $U_j\subseteq U_{j+1}$ and $V_j\subseteq V_{j+1}$. Suppose that there is a surjective map $T\colon U\to V$ such that
\[
  T|_{U_j}\colon U_j\to V_j
\]
is surjective for all $j\geq 0$ and constants $M,\alpha>0$, independent of $j$, such that
\begin{align}
  \norm{Tx-Ty}_V
  &\le M\norm{x-y}_U,
    \label{eq:pg-map-bound}\\
  \dual{Ax-Ay}{Tx-Ty}
  &\ge\alpha\norm{x-y}_U^2,
    \qquad x,y\in U_j.
    \label{eq:pg-monotonicity}
\end{align}
We consider the solutions $u_j\in U_j$ to
\begin{align}\label{eq:discrete}
    \dual{Au_j}{v}=\dual{f}{v}\quad\text{for all }v\in V_j.
\end{align}

The well-posedness follows standard arguments and is included for completeness only.
\begin{lemma}[Discrete well-posedness]\label{lem:pg-wellposed}
For every $f\in V^\star$, the equations~\eqref{eq:discrete}
have a unique solution.
\end{lemma}

\begin{proof}
For $x,y\in U_j$, the assumptions give
\begin{align}\label{eq:invLipschitz}
  \alpha\norm{x-y}_U^2
  \le\norm{Ax-Ay}_{V^\star}\norm{Tx-Ty}_V
  \le L\norm{x-y}_U\norm{Tx-Ty}_V.
\end{align}
This proves that $T_j:=T|_{U_j}$ is injective. By its assumed surjectivity it is
a bijection, and its inverse is $(L/\alpha)$-Lipschitz.
If $V_j=\{0\}$, bijectivity implies $U_j=\{0\}$ and the assertion
is immediate. Otherwise, define $B_j:V_j\to V_j^\star$ by
\[
  \dual{B_jz}{v}:=\dual{A(T_j^{-1}z)}{v},
  \qquad z,v\in V_j.
\]
Writing $x=T_j^{-1}z$ and $y=T_j^{-1}\zeta$, we obtain
\begin{align*}
  \norm{B_jz-B_j\zeta}_{V_j^\star}
  &\le\norm{Ax-Ay}_{V^\star}
   \le\frac{L^2}{\alpha}\norm{z-\zeta}_V,\\
  \dual{B_jz-B_j\zeta}{z-\zeta}
  &\ge\alpha\norm{x-y}_U^2
   \ge\frac{\alpha}{M^2}\norm{z-\zeta}_V^2.
\end{align*}
Set $\Lambda=L^2/\alpha$ and $\beta=\alpha/M^2$.
These inequalities imply $0<\beta\le\Lambda$. Let
$R_j:V_j\to V_j^\star$ be the Riesz isometry, and put $f_j=f|_{V_j}$.
The map
\[
  F_j(z):=z-\frac\beta{\Lambda^2}R_j^{-1}(B_jz-f_j),
  \qquad z\in V_j
\]
satisfies
\begin{align*}
  \norm{F_j(z)-F_j(\zeta)}_V^2
  &=\norm{z-\zeta}_V^2
    -\frac{2\beta}{\Lambda^2}
       \dual{B_jz-B_j\zeta}{z-\zeta}
    +\frac{\beta^2}{\Lambda^4}
       \norm{B_jz-B_j\zeta}_{V_j^\star}^2\\
  &\le\left(1-\frac{\beta^2}{\Lambda^2}\right)
       \norm{z-\zeta}_V^2.
\end{align*}
The space $V_j$ is complete, so the contraction principle gives a
unique fixed point $z$, equivalently a unique solution of $B_jz=f_j$.
Then $u_j=T_j^{-1}z$ is the unique solution.
\end{proof}
\noindent
We show the assumptions~\eqref{eq:cea} and~\eqref{eq:sep} for the sequence $\{Tu_1,\ldots, Tu_n\}\subseteq V$.

\textbf{C\'ea type estimate:}
We obtain for all $0\leq i\leq j\leq k\leq \ell\leq n$
\begin{align*}
    \alpha \norm{u_k-u_j}_U^2\leq \dual{Au_k-Au_j}{Tu_k-Tu_j} = \dual{Au_k-Au_j}{Tu_k-Tu_i}\leq LM\norm{u_k-u_j}_U\norm{u_k-u_i}_U
\end{align*}
and 
\begin{align*}
    \alpha \norm{u_k-u_i}_U^2\leq \dual{Au_k-Au_i}{Tu_k-Tu_i} = \dual{Au_\ell-Au_i}{Tu_k-Tu_i}\leq LM\norm{u_\ell-u_i}_U\norm{u_k-u_i}_U.
\end{align*}
Both estimates can be combined with the bi-Lipschitz property of $T$ to show~\eqref{eq:cea} with $C_{\text{C\'ea}}=(LM/\alpha)^3$.

\textbf{Separation:} To see~\eqref{eq:sep}, consider $v\in {\rm conv}\{Tu_i\,|\,1\leq i\leq n,\,i\leq J\}$ and
    \begin{align*}
        \dual{Au_{\max J}-Au_{\min J}}{\mu_{J,\lambda}-v} &= 
         \dual{Au_{\max J}-Au_{\min J}}{\mu_{J,\lambda}-Tu_{\min J}} \\
         &= \sum_{j\in J}\lambda_{j}\dual{Au_j-Au_{\min J}}{Tu_j-Tu_{\min J}} \\
         &\geq 
         \alpha/M^2\sum_{j\in J} \lambda_{j}\norm{Tu_j-Tu_{\min J}}_V^2\geq \alpha/M^2 E_J,
         \end{align*}
    where the last inequality follows since the mean minimizes the squared distances. Cauchy--Schwarz yields
    \begin{align*}
     L\norm{u_{\max J}-u_{\min J}}_{U} \norm{\mu_{J,\lambda}-v}_V \geq \norm{Au_{\max J}-Au_{\min J}}_{V^\star}\norm{\mu_{J,\lambda}-v}_V\geq \alpha/M^2 E_{J}.
    \end{align*}
    Moreover, we have for $i=\min J$ and $j=\max J$ that
    \begin{align*}
        E_{ J}\geq \frac12\Big(\norm{Tu_i-(Tu_i+Tu_j)/2}_V^2+\norm{Tu_j-(Tu_i+Tu_j)/2}_V^2\Big)\geq\alpha^2/(4L^2) \norm{u_{\max J}-u_{\min J}}_{U}^2,
    \end{align*}
    where the second inequality follows by \eqref{eq:invLipschitz}
    and hence separation~\eqref{eq:sep} with $C_{\rm sep}:=4L^4M^4/\alpha^4$.
  Thus Theorem~\ref{thm:qo} applies to this sequence and by the bi-Lipschitz property of $T$ also yields the result for the sequence $\{u_1,\ldots, u_n\}\subseteq U$.

\section{Application to nonlinear ODEs}\label{sec:examples}
\noindent
In the following we introduce continuous Petrov--Galerkin time stepping for ordinary differential equations. This example showcases, in which applications the abstract Theorem \ref{thm:qo} is particularly useful. 

\noindent For a bounded interval $[t_0,t_{\rm end}] \subset \R$, we consider the initial value problem
\begin{align*}
    y'(t)&=F(t,y(t))\quad\text{for every }t\in [t_0,t_{\rm end}],\\
     y(t_0)&=x_0.
\end{align*}
The right-hand side $F:[t_0,t_{\rm end}]\times\R^d\to\R^d$ is supposed to be continuous in its first
argument and Lipschitz continuous in its second, i.e. for some $L_F>0$,
\begin{equation}\label{eq:ode-weighted-lipschitz}
  |F(t,a)-F(t,b)|\le L_F|a-b|,
  \qquad a,b\in\R^d, \quad t \in [t_0,t_{\rm end}].
\end{equation}

\subsection{Continuous Petrov--Galerkin time stepping}
\noindent
Let $(\mathcal I_j)_{j\ge0}$ be a sequence of partitions of $[t_0,t_{\rm end}]$ into finitely many intervals, where each partition $\mathcal I_{j+1}$ is obtained from $\mathcal I_j$ by subdividing some of its intervals. Such sequences are generated, e.g., by adaptive step-size algorithms. We denote the nodes of the initial partition $\mathcal I_0$ by $t_0=t_0^0<t_1^0<\cdots<t_m^0=t_{\rm end}$ and its maximal step size by
\[
  \tau_0:=\max_{1\le n\le m}(t_n^0-t_{n-1}^0).
\]
We require the initial partition to be fine enough in the sense that
\begin{equation}\label{eq:ode-step}
  L_F\tau_0\le\frac{\log 2}{4}.
\end{equation}

\noindent
We work with the spaces
\[
  U:=\{w\in H^1(t_0,t_{\rm end};\R^d):w(t_0)=0\},
  \qquad V:=L^2(t_0,t_{\rm end};\R^d),
\]
whose norms are specified below. For a fixed polynomial degree $p\in\N$, the trial and test spaces are
\begin{align*}
  U_j&:=\{w\in U:w|_I\in\mathcal P^p(I;\R^d)
                         \text{ for every }I\in\mathcal I_j\},\\
  V_j&:=\{v\in V:v|_I\in\mathcal P^{p-1}(I;\R^d)
                         \text{ for every }I\in\mathcal I_j\}.
\end{align*}
Hence, the trial functions are continuous piecewise polynomials of degree $p$, whereas the test functions are piecewise polynomials of degree $p-1$ which are allowed to jump at the nodes. Both sequences of spaces are nested and finite dimensional. The continuous Petrov--Galerkin method seeks an approximation $y_j$ of $y$ such that
\begin{equation}\label{eq:ode-cg}
  y_j-x_0\in U_j,\qquad
  \int_{t_0}^{t_{\rm end}}\bigl(y_j'(t)-F(t,y_j(t))\bigr)\cdot v_j(t)\,dt=0
  \quad\text{for every }v_j\in V_j.
\end{equation}

\subsection{Application of the abstract result}
\noindent
In order to apply Theorem~\ref{thm:qo}, we have to choose the operator $A$, the norms on $U$ and $V$, and the map $T$. We incorporate the initial value in the right-hand side and define $A:U\to V^\star$ by
\begin{equation}\label{eq:ode-operator}
  \dual{Au}{v}
  :=\int_{t_0}^{t_{\rm end}}\bigl(u'(t)-F(t,x_0+u(t))\bigr)\cdot v(t)\,dt.
\end{equation}
The assumptions on $F$ imply that $t\mapsto F(t,x_0+u(t))$ is continuous and bounded for every $u\in U$, so $A$ is well defined. By the Picard--Lindel\"of theorem, the initial value problem has a unique solution $y\in C^1([t_0,t_{\rm end}];\R^d)$. Thus, $u:=y-x_0\in U$ satisfies $Au=0$, and $u_j:=y_j-x_0\in U_j$ solves the Petrov--Galerkin equations~\eqref{eq:discrete} with $f=0$ if and only if $y_j$ solves~\eqref{eq:ode-cg}.

\noindent
The obvious choice of the remaining ingredients is to equip $U$ with $\norm{w'}_{L^2(t_0,t_{\rm end})}$, to equip $V$ with the $L^2$-norm, and to set $Tw:=w'$. This satisfies $T(U_j)=V_j$. However, for $x,z\in U$ and $w:=x-z$, the Lipschitz bound~\eqref{eq:ode-weighted-lipschitz} and the Poincar\'e inequality only show
\[
  \dual{Ax-Az}{w'}
  \ge\norm{w'}_{L^2(t_0,t_{\rm end})}^2-L_F\norm{w}_{L^2(t_0,t_{\rm end})}\norm{w'}_{L^2(t_0,t_{\rm end})}
  \ge\Bigl(1-\frac{L_F(t_{\rm end}-t_0)}{\sqrt2}\Bigr)\norm{w'}_{L^2(t_0,t_{\rm end})}^2.
\]
This yields the strong monotonicity~\eqref{eq:pg-monotonicity} only for short time intervals with $L_F(t_{\rm end}-t_0)<\sqrt2$. To avoid this restriction, we use exponentially decaying weights, first introduced by \cite{bielecki}, as is common for Gronwall-type arguments.

\noindent
We define the piecewise constant weight
\begin{equation}\label{eq:ode-weight}
  \omega(t):=e^{-8L_F(t_{n-1}^0-t_0)}
  \qquad\text{for }t_{n-1}^0<t<t_n^0,\quad n=1,\ldots,m.
\end{equation}
The weight is fixed on the initial partition $\mathcal I_0$ and does not depend on the level $j$. We equip $U$ and $V$ with the norms
\begin{equation}\label{eq:ode-norms}
  \norm{w}_U^2:=\int_{t_0}^{t_{\rm end}}\omega(t)|w'(t)|^2\,dt,
  \qquad
  \norm{v}_V^2:=\int_{t_0}^{t_{\rm end}}\omega(t)^{-1}|v(t)|^2\,dt.
\end{equation}
Since $e^{-8L_F(t_{\rm end}-t_0)}\le\omega\le1$, these norms are equivalent to the standard norms and $U$, $V$ are Hilbert spaces. Moreover, for $r\in L^2(t_0,t_{\rm end};\R^d)$, the functional $v\mapsto\int_{t_0}^{t_{\rm end}}r\cdot v\,dt$ satisfies
\begin{equation}\label{eq:ode-dual-norm}
  \norm{r}_{V^\star}^2=\int_{t_0}^{t_{\rm end}}\omega(t)|r(t)|^2\,dt.
\end{equation}

\noindent
The gain of the weighted norms is a Poincar\'e inequality whose constant does not depend on the length of the time interval. To see this, let $\rho(t):=e^{-8L_F(t-t_0)}$. For $w\in U$, integration by parts and $w(t_0)=0$ give
\[
  \int_{t_0}^{t_{\rm end}}\rho\,w\cdot w'\,dt
  =\frac12\rho(t_{\rm end})|w(t_{\rm end})|^2
     +4L_F\int_{t_0}^{t_{\rm end}}\rho|w|^2\,dt
  \ge4L_F\int_{t_0}^{t_{\rm end}}\rho|w|^2\,dt.
\]
The Cauchy--Schwarz inequality therefore implies
\[
  \Bigl(\int_{t_0}^{t_{\rm end}}\rho|w|^2\,dt\Bigr)^{1/2}
  \le\frac{1}{4L_F}\Bigl(\int_{t_0}^{t_{\rm end}}\rho|w'|^2\,dt\Bigr)^{1/2}.
\]
On every interval of $\mathcal I_0$, there holds $\rho\le\omega\le e^{8L_F\tau_0}\rho$. Together with $e^{4L_F\tau_0}\le2$, which is the step-size restriction~\eqref{eq:ode-step}, we obtain
\begin{equation}\label{eq:ode-poincare}
  \Bigl(\int_{t_0}^{t_{\rm end}}\omega|w|^2\,dt\Bigr)^{1/2}
  \le\frac{e^{4L_F\tau_0}}{4L_F}
          \Bigl(\int_{t_0}^{t_{\rm end}}\omega|w'|^2\,dt\Bigr)^{1/2}
  \le\frac{1}{2L_F}\norm{w}_U
  \qquad\text{for every }w\in U.
\end{equation}

\noindent
We now verify the operator assumptions of Section~\ref{sec:pg}. Let $x,z\in U$ and $w:=x-z$. The identity~\eqref{eq:ode-dual-norm}, the Lipschitz bound~\eqref{eq:ode-weighted-lipschitz}, and the weighted Poincar\'e inequality~\eqref{eq:ode-poincare} show
\[
  \norm{Ax-Az}_{V^\star}
  \le\norm{w}_U+L_F\Bigl(\int_{t_0}^{t_{\rm end}}\omega|w|^2\,dt\Bigr)^{1/2}
  \le\frac32\norm{x-z}_U,
\]
which is the Lipschitz continuity~\eqref{eq:pg-lipschitz} with $L=3/2$. We choose
\[
  T:U\to V,\qquad Tw:=\omega w'.
\]
The definition of the norms~\eqref{eq:ode-norms} immediately gives $\norm{Tx-Tz}_V=\norm{x-z}_U$, which is~\eqref{eq:pg-map-bound} with $M=1$. Furthermore, every partition $\mathcal I_j$ is a refinement of $\mathcal I_0$, so $\omega$ is constant on each interval of $\mathcal I_j$. Therefore, $\omega w'\in V_j$ for all $w\in U_j$, and every $v\in V_j$ is the image of $w(t):=\int_{t_0}^t\omega^{-1}v\,ds\in U_j$. This shows surjectivity
\[
  T(U_j)=V_j\qquad\text{for every }j\ge0,
\]
 Finally, \eqref{eq:ode-weighted-lipschitz} and~\eqref{eq:ode-poincare} show
\begin{align*}
  \dual{Ax-Az}{Tx-Tz}
  &=\int_{t_0}^{t_{\rm end}}\omega\bigl(w'-F(t,x_0+x)+F(t,x_0+z)\bigr)\cdot w'\,dt\\
  &\ge\norm{w}_U^2
       -L_F\Bigl(\int_{t_0}^{t_{\rm end}}\omega|w|^2\,dt\Bigr)^{1/2}\norm{w}_U\\
  &\ge\frac12\norm{x-z}_U^2,
\end{align*}
which is the strong monotonicity~\eqref{eq:pg-monotonicity} with $\alpha=1/2$.

\noindent
Altogether, the assumptions~\eqref{eq:pg-lipschitz}--\eqref{eq:pg-monotonicity} are satisfied with
\[
  L=\frac32,\qquad M=1,\qquad
  \alpha=\frac12.
\]
Lemma~\ref{lem:pg-wellposed} with $f=0$ gives a unique solution $u_j\in U_j$ of~\eqref{eq:discrete} on each level, so $y_j=x_0+u_j$ is the unique solution of~\eqref{eq:ode-cg}.
Fix $\ell\in\Nzero$ and $n\in\N$. By Section~\ref{sec:pg}, the sequence $Tu_\ell,\ldots,Tu_{\ell+n}$ satisfies~\eqref{eq:cea} and~\eqref{eq:sep}.
Since $T$ is an isometry, $\norm{Tu_j-Tu_k}_V=\norm{u_j-u_k}_U=\norm{y_j-y_k}_U$ for all $j,k$.
Applying Theorem~\ref{thm:qo} to these $n+1$ points we obtain
\[
  \sum_{j=\ell}^{\ell+n-1}\norm{y_{j+1}-y_j}_U^2
  \lesssim n^{1-\delta}\norm{y_{\ell+n}-y_\ell}_U^2
  \qquad\text{with}\qquad
  \delta \in (0,1).
\]
It remains to return to the standard norm. Since $w(t_0)=0$ implies $\norm{w}_{L^2(t_0,t_{\rm end})}^2\le\frac12(t_{\rm end}-t_0)^2\norm{w'}_{L^2(t_0,t_{\rm end})}^2$, the bounds on the weight show
\[
  \norm{w}_U^2\le\norm{w}_{H^1(t_0,t_{\rm end})}^2
  \le\Bigl(1+\frac{(t_{\rm end}-t_0)^2}{2}\Bigr)e^{8L_F(t_{\rm end}-t_0)}\norm{w}_U^2
  \qquad\text{for every }w\in U.
\]
We conclude that, for every $\ell\in\Nzero$ and $n\in\N$,
\begin{equation}\label{eq:ode-qo}
  \sum_{j=\ell}^{\ell+n-1}\norm{y_{j+1}-y_j}_{H^1(t_0,t_{\rm end})}^2
  \lesssim n^{1-\delta}\norm{y_{\ell+n}-y_\ell}_{H^1(t_0,t_{\rm end})}^2.
\end{equation}
The implied constant and $\delta$ as well as the step-size restriction~\eqref{eq:ode-step} are independent of the polynomial degree $p$, the dimension $d$, and the partitions $\mathcal I_j$.

\subsection{Comparison with existing results}
\noindent
A quasi-orthogonality estimate for the method~\eqref{eq:ode-cg} was previously shown in~\cite{adaptive_ode}. That work considers an adaptive algorithm which is steered by a residual error estimator $\eta_j$ and bisects the time steps selected by D\"orfler marking.  As a key ingredient in their analysis, \cite[Lemma~12]{adaptive_ode} proves the following: Suppose that, in addition to~\eqref{eq:ode-weighted-lipschitz}, the right-hand side $F$ is continuously differentiable with Lipschitz continuous Jacobian, and that the initial partition $\mathcal I_0$ is sufficiently fine. Then, there exist $\varepsilon>0$, $1>\delta>0$, and $C(n)=\mathcal O(n^{1-\delta})$ such that
\[
  \sum_{j=\ell}^{\ell+n-1}\Bigl(\norm{y_{j+1}-y_j}_{H^1(t_0,t_{\rm end})}^2-\varepsilon\,\eta_j^2\Bigr)
  \le C(n)\norm{y-y_\ell}_{H^1(t_0,t_{\rm end})}^2
  \qquad\text{for all }\ell,n\in\N.
\]
The proof linearizes the initial value problem around the exact solution $y$ and shows that the continuous Petrov--Galerkin discretization of the linearized problem is uniformly inf-sup stable on sufficiently fine partitions. Thus, the linear result of~\cite{infsup} applies to the linearized problem. The difference between the discrete solutions of the nonlinear and the linearized problem is quadratic in the error $\norm{y-y_j}_{H^1(t_0,t_{\rm end})}$ and is absorbed by the term $\varepsilon\,\eta_j^2$.

\noindent
The estimate~\eqref{eq:ode-qo} differs from this result in the following ways. First, we apply Theorem~\ref{thm:qo} directly to the nonlinear operator~\eqref{eq:ode-operator}. No linearization is necessary, and the Lipschitz continuity~\eqref{eq:ode-weighted-lipschitz} is the only assumption on $F$, whereas the Taylor expansion in~\cite{adaptive_ode} requires the differentiability of $F$. Second, \eqref{eq:ode-qo} does not contain an error estimator and holds for arbitrary nested refinements of an initial partition satisfying \eqref{eq:ode-step}. In particular, it implies the estimate of~\cite[Lemma~12]{adaptive_ode} with $\varepsilon=0$. Third, the step-size restriction~\eqref{eq:ode-step} is explicit and depends only on the Lipschitz constant $L_F$.
\printbibliography

\end{document}